\documentclass[12pt,reqno]{amsart}

\usepackage{amssymb,latexsym}

\usepackage{enumerate}
\allowdisplaybreaks
\usepackage[french,english]{babel}
\usepackage{amsmath}
\usepackage{graphicx}
\usepackage{amssymb}
\usepackage{bbm}
\usepackage{amsthm,mathtools}
\usepackage{ulem}
\usepackage{geometry}
\usepackage{tikz-cd}
\usepackage{mathrsfs}
\usepackage[colorinlistoftodos]{todonotes}
\usepackage{enumitem}
\usepackage{verbatim}
\usepackage[foot]{amsaddr}
\usepackage{dsfont}
\usepackage{cite}
\usepackage[T1]{fontenc}

\makeatletter

\@namedef{subjclassname@2010}{
	
	\textup{2020} Mathematics Subject Classification}

\makeatother
\newtheorem{thm}{Theorem}[section]
\newtheorem*{thm*}{Theorem}

\theoremstyle{definition}

\newtheorem{Remark}{Remark}[section]
\numberwithin{equation}{section}

\newcommand{\Repart}{\operatorname{Re}} 

\newcommand{\inv}{^{-1}}

\newcommand{\mbz}{\mathbb{Z}}
\newcommand{\mbr}{\mathbb{R}}

\newcommand{\mbn}{\mathbb{N}}

\newcommand{\mce}{\mathcal{E}}

\usepackage{hyperref}
\hypersetup{colorlinks=true,linkcolor=blue,anchorcolor=blue,citecolor=blue}
\usepackage{color}
\newcommand{\newabstract}[1]{%
	\par\bigskip
	\csname otherlanguage*\endcsname{#1}%
	\csname captions#1\endcsname
	\item[\hskip\labelsep\scshape\abstractname.]
}

\begin{document}

	\baselineskip=17pt

	\title[$\Omega$-results for the logarithmic derivative of $\zeta(s)$ on vertical homogeneous progressions]{$\Omega$-results for the  logarithmic derivative of $\zeta(s)$ on vertical homogeneous progressions}

    \author{Shengbo Zhao\textsuperscript{1}}
	\address{1.School of Mathematical Sciences, Key Laboratory of Intelligent Computing and Applications(Ministry of Education), Tongji University, Shanghai 200092, P. R. China}	\email{shengbozhao@hotmail.com}

	\begin{abstract} 
	In this paper, we establish $\Omega$-results for the logarithmic derivative of the Riemann zeta function on vertical homogeneous progressions on the 1-line and in the critical strip. Compared with Yang's work in 2023, our results show that the discrete case has a similar order of magnitude to the continuous case.
	\end{abstract}
	
    \keywords{$\Omega$-results, the Riemann zeta function, the logarithmic derivative, homogeneous progressions, resonance method. }
	
	\subjclass[2020]{Primary 11M06, 11N37.}
	
	\maketitle

\section{Introduction}
The Riemann zeta function $\zeta(s)$ has always been one of the most important research objects in the analytic number theory. Its logarithmic derivative, $\zeta^\prime /\zeta(s)$, also plays a fundamental role in numerous studies. In particular,it is closely related to the Chebyshev function $\psi(x) = \sum_{n \le x}\Lambda(n)$, which serves as a key tool in understanding the distribution of prime numbers. Moreover, $\zeta^\prime /\zeta(s)$ appears in the investigation of zero-density results for $\zeta(s)$.
\par
Here and throughout, we let $s=\sigma +it$ denote a complex variable, and write $\log_j$ for the $j$-th iterated logarithm. In 2023, Yang \cite{yang2023omega} established $\Omega$-results for $\zeta^\prime/\zeta(s)$ for $\sigma \in (1/2,1]$ and $t \in [T,T^\beta]$, where $\beta \in (0,1)$ is fixed. For sufficiently large $T$ and any fixed $\varepsilon \in (0,1)$, he showed that
$$
\max_{T^\beta \le t\le T} -\Repart\frac{\zeta^\prime}{\zeta}(1+it) \ge \log_2 T+\log_3 T + C -\varepsilon,
$$
and for $\sigma \in (1/2,1)$,
$$
 \max_{T^\beta \le t \le T} -\Repart\frac{\zeta^\prime}{\zeta}(\sigma+it) \ge C(\sigma)(\log T)^{1-\sigma}(\log_2 T)^{1-\sigma},
$$
where $C$ and $C(\sigma)$ are some constant. Notably, his celebrated work improved bounds of Landau \cite{Landau1911}, and Bohr and Landau \cite{Bohr1913}.
\par
Investigating the behavior of $\zeta(s)$ in vertical arithmetic progressions has attracted significant attention due to its applications in analytic number theory and probabilistic models. Li and Radziwiłł \cite{li2015theriemann} provided that the moments, zero distributions, and various other properties of $\zeta(s)$ over vertical arithmetic progressions are closely related to those in the continuous case on the critical line. Furthermore, Minelli and Sourmelidis \cite{minelli2025discrete} showed that large values of $\zeta(s)$ can be captured along any given homogeneous progression, with bounds of the same quality as the best currently known $\Omega$-results in the continuous case. 
\par
Motivated by the work of \cite{minelli2025discrete,yang2023omega}, we give $\Omega$-results for $\zeta^\prime/\zeta(s)$ when $\sigma \in (1/2,1]$. Our results demonstrate that $\Omega$-results on homogeneous progressions achieve the same order of magnitude as those in the continuous case. Henceforth, we assume that $\ell \in \mbz$, where $\mbz$ is the set of all integers. First, we give the following $\Omega$-results on the 1-line.
\begin{thm}
    \label{thm1}
    Let $\alpha>0$, $\theta,\varepsilon \in (0,1)$ be fixed. Let $N$ be sufficiently large. Then we have
    \begin{equation*}
        \max_{N^\theta \le \ell \le N} -\Repart\frac{\zeta^\prime}{\zeta}(1+i\alpha\ell) \ge \log_2 N+\log_3 N +c(\theta) -\varepsilon,
    \end{equation*}
    where $c(\theta)$ is a negative constant depending on $\theta$. Specifically,
    \begin{equation*}
        c(\theta) = \log(1-\theta) - \log_2 4 -\gamma -\sum_{k=2}^{\infty}\sum_p\frac{\log p}{p^k}-1.
    \end{equation*}
\end{thm}

\begin{Remark}
    Let $A$ be any positive real number. Motivated by \cite{yang2022logtype}, we can get $\Omega$-results for $\zeta^\prime/\zeta(\sigma_A +i\alpha\ell)$, where $\sigma_A = 1 - A/\log_2 N$. This result introduces a coefficient $m(A)$ in front of the term $\log_2 N$, with $m(A) \to 1$ as $A \to 0$, .
\end{Remark}

\begin{Remark}
    Following the proof of \cite[Theorem 1.2]{yang2023omega}, both for $\sigma=1$ and $\sigma \in (1/2,1)$ we can derive measure estimates for the set on which these large values occur. Since the details are nearly identical, we omit these results.
\end{Remark}

We next state $\Omega$-results in the critical strip. 
\begin{thm}
    \label{thm2}
    Let $\alpha>0$, $\theta \in (0,1)$ and $\sigma \in (1/2,1)$ be fixed. Let $N$ be sufficiently large. Then we have
    \begin{equation*}
        \max_{N^\theta \le \ell \le N} -\Repart\frac{\zeta^\prime}{\zeta}(\sigma+i\alpha\ell) \ge C(\sigma, \theta)(\log N)^{1-\sigma}(\log_2 N)^{1-\sigma},
    \end{equation*}
    where $C(\sigma,\theta)$ is some positive constant which can be effectively computed.
\end{thm}

\begin{Remark}
    For any fixed $\delta \in [0, 2\pi]$ and $\sigma \in (1/2, 1]$, we can also establish $\Omega$-results for $\Repart(e^{-i\delta}\zeta^\prime/\zeta(\sigma+i\alpha\ell))$. As expected, these lower bounds are not as sharp as those in the case of $\delta = \pi$, that is, our Theorems \ref{thm1} and \ref{thm2}. Therefore, we omit these weaker results in this paper.
\end{Remark}

Our proofs rely on the resonance method. The central idea is to construct a suitable Dirichlet series $R(t)$ as a resonator, which resonates with the target function and picks out its large values. For further discussion of the resonance method and related results, we recommend \cite{aistleitner2016lower,aistleitner2019extreme,aistleitner2019onlarge,bondarenko2017large,bondarenko2018argument,soundararajan2008extreme,qiyu2024large} and the references therein.
\par
Finally, we introduce some notations. Let $p$ be a prime number. Let $\varepsilon > 0$ be an arbitrarily small number. We remark that each occurrence of $\varepsilon$ may represent a different value. Furthermore, we denote the Fourier transform of a function $f \in L^1(\mathbb{R})$ as
$$\widehat{f}(\xi) \coloneqq \int_{\mathbb R} f(x) e^{-2\pi i  \xi x} \mathrm{d}x.$$ 

\section{Proof of Theorem \ref{thm1}}
\label{section2}
Set $K = \exp((\log N)^2)$ and $\Phi(y) \coloneqq e^{-y^{2}/2}$ as in \cite{bondarenko2017large}. Plainly, the Fourier transform $\widehat{\Phi}$ satisfies $\widehat{\Phi}(\xi) = \sqrt{2\pi}\Phi(2\pi\xi) >0$ for all $\xi \in \mbr$. In the subsequent proof, we will frequently make use of the positivity of $\Phi$ and $\widehat{\Phi}$. Define the sums
\begin{align*}
    S_1 &\, \coloneqq S_1(R,N)= \sum_{\ell \in \mbz} |R(\alpha \ell)|^2 \Phi \Big(\frac{\ell \log N}{N} \Big), \\
    S_2 &\, \coloneqq S_2(R,N)= \sum_{\ell \in \mbz} \Repart\Big(\sum_{n\le K}\frac{\Lambda(n)}{n^{1+i\alpha\ell}} \Big) |R(\alpha \ell)|^2 \Phi \Big(\frac{\ell \log N}{N} \Big).
\end{align*}
Let $M = \kappa\log N \log_2 N$, where $\kappa = (1-\theta)e^{-\varepsilon}/(\log 4)$. Then as in \cite{aistleitner2019onlarge}, define a completely multiplicative function $r(n)$ by $r(p) = 1-p/M$, when $p \le M$. Furthermore, if $p>M, \,r(p)=0$. Define the resonator $R(t) \coloneqq \prod_p (1-r(p)p^{-it})\inv=\sum_{n \in \mbn}r(n)n^{-it}$. The prime number theory leads to 
\begin{equation*}
    \log|R(t)| \le \log M\sum_{p \le M}1 - \sum_{p \le M}\log p \le \frac{M}{\log M}+ O\Big( \frac{M}{(\log M)^2}\Big). 
\end{equation*}
Combining with $M = \kappa\log N \log_2 N$, we conclude that
\begin{equation}
    \label{Rupperbound1}
    |R(t)|^2 \le N^{2\kappa+o(1)}.
\end{equation}
\par
Noting that $\widehat{\Phi}$ is always positive, we can obtain the following lower bound for $S_1$
\begin{equation*}
     S_1  = \frac{N}{\log N}\sum_{m,n \in \mbn}r(m)r(n)\sum_{\ell \in \mbz} \widehat{\Phi}\Big(\frac{N}{\log N}\Big(\frac{\alpha}{2\pi}\log \frac{n}{m}-\ell \Big) \Big) \ge \sqrt{2\pi}\frac{N}{\log N}\sum_{n \in \mbn}r(n)^2 
\end{equation*}
by applying the Poisson summation formula. Using the results for \cite[p. 841]{aistleitner2019onlarge}, we have
\begin{equation}
    \label{S1lowerbound}
    S_1 \ge N^{1+(2-\log 4)\kappa+o(1)}.
\end{equation}
\par
On the other hand, for $S_2$, 
\begin{align*}
    S_2 & \, = \sum_{p^k \le K}\frac{\log p}{p^k}\sum_{m,n \in \mbn}r(m)r(n)\sum_{\ell \in \mbz}\Big( \frac{n}{p^km}\Big)^{i\alpha \ell}\Phi\Big(\frac{\ell \log N}{N} \Big) \\
    & \, \ge \sum_{p \le M}\frac{\log p}{p}\sum_{m,n \in \mbn}r(m)r(n)\sum_{\ell \in \mbz}\Big( \frac{n}{pm}\Big)^{i\alpha \ell}\Phi\Big(\frac{\ell \log N}{N} \Big) \\
    & \, \ge \sum_{p \le M}\frac{\log p}{p}r(p)\sum_{m,n \in \mbn}r(m)r(n)\sum_{\ell \in \mbz}\Big( \frac{n}{m}\Big)^{i\alpha \ell}\Phi\Big(\frac{\ell \log N}{N} \Big).
\end{align*}
Hence, it is clear that
\begin{equation}
    \label{S2S1ratio}
    \frac{S_2}{S_1} \ge \sum_{p \le M}\frac{\log p}{p}r(p).
\end{equation}
\par
By \cite[Eq. (1)]{yang2023omega}, uniformly for $t \in [\alpha N^\theta, \alpha N]$, we have 
\begin{equation}
    \label{appro1}
    \sum_{n \le K}\frac{\Lambda(n)}{n^{1+it}} = - \frac{\zeta^\prime}{\zeta}(1+it) + O(N^{-B(\theta)}).
\end{equation}
Here, $B(\theta)$ is a positive constant, which depends on $\theta \in (0,1)$. Trivially, 
\begin{equation}
    \label{Reupper1}
    \Big|\Repart\Big(\sum_{n\le K}\frac{\Lambda(n)}{n^{1+it}} \Big) \Big| \le \sum_{n \le K}\frac{\Lambda(n)}{n}\ll K \ll N^{o(1)}.
\end{equation}
Thus, combining with \eqref{Rupperbound1}, we have
$$
\Big|\sum_{|\ell|\le N^\theta}\Repart\Big(\sum_{n\le K}\frac{\Lambda(n)}{n^{1+i\alpha \ell}} \Big)|R(\alpha \ell)|^2\Phi\Big(\frac{\ell \log N}{N}\Big) \Big| \le N^{\theta+2\kappa+o(1)}.
$$
Furthermore, \eqref{Rupperbound1}, \eqref{Reupper1} and the rapid decay of $\Phi$ imply that
$$
\Big|\sum_{|\ell|\ge N}\Repart\Big(\sum_{n\le K}\frac{\Lambda(n)}{n^{1+i\alpha \ell}} \Big)|R(\alpha \ell)|^2\Phi\Big(\frac{\ell \log N}{N}\Big) \Big| \le 1.
$$
Similarly, 
$$
\Big|\sum_{|\ell|\le N^\theta}|R(\alpha \ell)|^2\Phi\Big(\frac{\ell \log N}{N}\Big) \Big| \le N^{\theta+2\kappa+o(1)} \quad \text{and} \quad \Big|\sum_{|\ell|\ge N}|R(\alpha \ell)|^2\Phi\Big(\frac{\ell \log N}{N}\Big) \Big| \le 1.
$$
At this point, \eqref{S1lowerbound} and \eqref{appro1} provide that
$$
\max_{N^\theta \le \ell \le N} -\Repart\frac{\zeta^\prime}{\zeta}(1+i\alpha\ell) \ge \frac{S_2}{S_1} + O(N^{\theta+(\log 4)\kappa-1+o(1)}).
$$
\par
By the definition of $r(n)$ and \cite[p. 68]{RosserSchoenfeld1962}, we directly deduce that
$$
\sum_{p \le M}\frac{\log p}{p}r(p) = \log M - \gamma -\sum_{k=2}^{\infty}\sum_p\frac{\log p}{p^k}-1 + O\Big(e^{-B\sqrt{\log M}} \Big),
$$
where $\gamma$ is the Euler-Mascheroni constant, and $B$ is some absolute positive constant. Hence, \eqref{S2S1ratio} yields
$$
\max_{N^\theta \le \ell \le N} -\Repart\frac{\zeta^\prime}{\zeta}(1+i\alpha\ell) \ge \log M - \gamma -\sum_{k=2}^{\infty}\sum_p\frac{\log p}{p^k}-1 + \mce,
$$
where $$\mce =  O\Big(e^{-B\sqrt{\log M}} \Big)+O(N^{\theta+(\log 4)\kappa-1+o(1)}).$$
Finally, substituting $M = \kappa\log N \log_2 N$ and $\kappa = (1-\theta)e^{-\varepsilon}/(\log 4)$, the proof of Theorem \ref{thm1} is completed.

\section{Proof of Theorem \ref{thm2}}
In contrast to Section \ref{section2}, we now set $K=(\log N)^{20/\varepsilon}$, where $\varepsilon \in (0,\sigma-1/2)$ is fixed.
Combining a classical estimate of Ingham \cite{Ingham1940} with \cite[Lemma 1]{yang2023omega}, we can derive the following effective asymptotic formula:
\begin{equation}
    \label{appro2}
    -\frac{\zeta^\prime}{\zeta}(\sigma+it) = \sum_{n \le K}\frac{\Lambda(n)}{n^{\sigma+it}}+O\Big((\log N)^{-18} \Big),\quad \forall t \in [\alpha N^\theta,\alpha N]\,\setminus\, \mce(\sigma,N),
\end{equation}
where the exceptional set $\mce(\sigma,N)$ has measure
$$
\operatorname{meas}(\mce(\sigma,N)) \ll N^{\frac{3(1-\sigma+\varepsilon)}{2-\sigma+\varepsilon}}(\log N)^{5+\frac{20}{\varepsilon}}.
$$
\par
Define $M=\kappa \log N \log_2 N$, where $\kappa>0$ will be chosen later. Then as in \cite{XiaoYang2022}, let $r(n)$ be a completely multiplicative function with values for $p\le M$ as $r(p) = 1-(p/M)^\sigma.$
Furthermore, if $p>M, \,r(p)=0$. Define the resonator $R(t) \coloneqq \prod_p (1-r(p)p^{-it})\inv=\sum_{n \in \mbn}r(n)n^{-it}$, too. Similarly, the prime number theorem gives
\begin{equation}
    \label{Rupperbound2}
    |R(t)|^2 \le N^{2\sigma\kappa+o(1)}.
\end{equation}
We consider the sums $S_1$, which is defined as in Section \ref{section2}, and
$$
S_2(\sigma) \coloneqq S_2(\sigma,R,N)= \sum_{\ell \in \mbz} \Repart\Big(\sum_{n\le K}\frac{\Lambda(n)}{n^{\sigma+i\alpha\ell}} \Big) |R(\alpha \ell)|^2 \Phi \Big(\frac{\ell \log N}{N} \Big),
$$
where $\Phi(y) \coloneqq e^{-y^{2}/2}$. Applying the Poisson summation formula and using the computation in \cite[p. 79]{Dong2022}, we find
\begin{equation}
    \label{S1lower}
    S_1 \ge \sqrt{2\pi}\frac{N}{\log N}\sum_{n \in \mbn}r(n)^2 \ge N^{1+(1-\lambda(\sigma))\sigma\kappa},
\end{equation}
where $\lambda(\sigma) \coloneqq \int_0^1 t^\sigma/(2-t^\sigma)\mathrm{d}t$ is a small positive quantity depending only on $\sigma$.
\par
As in the proof of Theorem \ref{thm1}, an analogous approach of the Dirichlet polynomial shows that
\begin{equation}
    \label{S2sigmaS1ratio}
    S_2(\sigma) \ge S_1\sum_{p \le M}\frac{\log p}{p^\sigma}r(p).
\end{equation}
However, unlike Section \ref{section2}, we must also control the contribution from the exceptional set $\mce(\sigma,N)$. Employing \eqref{Rupperbound2}, the contribution is at most $O(N^{2\sigma\kappa+\frac{3(1-\sigma+\varepsilon)}{2-\sigma+\varepsilon}+o(1)}).$
To ensure that this term is negligible compared $S_1$, we impose the requirement that
\begin{equation}
    \label{kappa1}
    2\sigma\kappa+\frac{3(1-\sigma+\varepsilon)}{2-\sigma+\varepsilon} < 1+(1-\lambda(\sigma))\sigma\kappa.
\end{equation}
Since $K=(\log N)^{20/\varepsilon}$, we crudely have
\begin{equation}
    \label{crude}
     \Big|\Repart\Big(\sum_{n\le K}\frac{\Lambda(n)}{n^{\sigma+it}} \Big) \Big| \le \sum_{n \le K}\frac{\Lambda(n)}{n^\sigma}\ll K \ll N^{o(1)}.
\end{equation}
Assume \eqref{kappa1}, we obtain from \eqref{appro2}, \eqref{S2sigmaS1ratio} and \eqref{crude} that
\begin{equation*}
     \max_{N^\theta \le \ell \le N} -\Repart\frac{\zeta^\prime}{\zeta}(\sigma+i\alpha\ell) \ge \sum_{p \le M}\frac{\log p}{p^\sigma}r(p)+ O(N^{\theta +2\sigma\kappa-1-(1-\lambda(\sigma))\sigma\kappa+o(1)}).
\end{equation*}
The error term is small provided that
\begin{equation}
    \label{kappa2}
    \theta+2\sigma\kappa<1+(1-\lambda(\sigma))\sigma\kappa.
\end{equation}
From the prime number theorem and the choice $M=\kappa\log N \log_2 N$,
$$
\sum_{p \le M}\frac{\log p}{p^\sigma}r(p) = \Big(\frac{\sigma}{1-\sigma}+o(1) \Big)\kappa^{-\sigma}(\log N)^{1-\sigma}(\log_2 N)^{1-\sigma}.
$$
Choosing $\kappa$ satisfying \eqref{kappa1} and \eqref{kappa2}, we finish the proof of Theorem \ref{thm2}.

\begin{Remark}
    In the proof of Theorem \ref{thm2}, we make use of Ingham's classical zero-density theorem in \cite{Ingham1940}. In fact, when $\sigma$ is far from the critical line, sharper zero-density theorems may be used to enlarge the range of $\kappa$.
\end{Remark}

\begin{Remark}
    In Theorems \ref{thm1} and \ref{thm2}, our method can only yield results over homogeneous progressions $\alpha\ell \, (\alpha>0, \ell \in \mbz)$. When considering $\alpha\ell + \beta\,(\beta \neq 0)$, due to the existence of term $(n/pm)^{-i\beta}$ with an undetermined sign, we cannot directly discard many terms by using the positivity of $\Phi$ and $\widehat{\Phi}$. Therefore, more refined methods are required for more general arithmetic progressions.
\end{Remark}

	\bibliographystyle{siam}
    \bibliography{reference}

@article{aistleitner2019extreme,
  title={Extreme values of the {R}iemann zeta function on the 1-line},
  author={Aistleitner, Christoph and Mahatab, Kamalakshya and Munsch, Marc},
  journal={Int. Math. Res. Not.},
  volume={IMRN 2019},
  number={22},
  pages={6924--6932},
  year={2019}
}

@article{aistleitner2016lower,
  author = {C. Aistleitner},
  title = {Lower bounds for the maximum of the {R}iemann zeta function along vertical lines},
  journal = {Math. Ann.},
  volume = {365},
  number = {1-2},
  pages = {473--496},
  year = {2016}
}

@article{aistleitner2019onlarge,
    author = {Aistleitner, Christoph and Mahatab, Kamalakshya and Munsch, Marc and Peyrot, Alexandre},
    title = {On large values of ${L}(\sigma,
             \chi)$},
    journal = {Quart. J. Math.},
    year = {2019},
    volume = {70},
    pages = {831--848},
    number = {3}
}

@article{Bohr1913,
  author = {H. Bohr and E. Landau},
  title = {Beiträge zur {T}heorie der {R}iemannschen {Z}etafunktion},
  journal = {Math. Ann.},
  volume = {74},
  number = {1},
  pages = {3--30},
  year = {1913}
}

@article{bondarenko2017large,
  author = {A. Bondarenko and K. Seip},
  title = {Large greatest common divisor sums and extreme values of the {R}iemann zeta function},
  journal = {Duke Math. J.},
  volume = {166},
  number = {9},
  pages = {1685--1701},
  year = {2017}
}

@article{bondarenko2018argument,
    author = {Bondarenko, A and Seip, K},
    title = {Extreme values of the {R}iemann zeta function and its argument},
    journal = {Math. Ann.},
    volume = {372},
    number = {3-4},
    pages = {999--1015},
    year = {2018}
}

@misc{Dong2022,
  author = {Z. Dong},
  title = {Distribution of values of the {R}iemann zeta function},
  howpublished = {Université Paris-Est Créteil Val-de-Marne - Paris 12},
  year = {2022},
  note = {Posted on 2022}
}

@article{Ingham1940,
  author = {A. E. Ingham},
  title = {On the estimation of ${N}(\sigma, T)$},
  journal = {Quart. J. Math.},
  volume = {11},
  year = {1940},
  pages = {291–292}
}

@article{Landau1911,
  author = {E. Landau},
  title = {Zur {T}heorie der {R}iemannschen {Z}etafunktion},
  journal = {Vierteljahrsschrift der Naturforschenden Gesellschaft in Zürich},
  volume = {56},
  year = {1911},
  pages = {125-148}
}

@article{RosserSchoenfeld1962,
  author = {J. B. Rosser and L. Schoenfeld},
  title = {Approximate formulas for some functions of prime numbers},
  journal = {Illinois J. Math.},
  volume = {6},
  pages = {64–94},
  year = {1962}
}

@article{li2015theriemann,
   author={Li, Xiannan and Radziwiłł, Maksym},
  journal={Int. Math. Res. Not. IMRN}, 
  title={The {R}iemann zeta function on vertical arithmetic progressions}, 
  year={2015},
  volume={2015},
  number={2},
  pages={325-354}
}

@article{minelli2025discrete,
    author = {Minelli, Paolo and Sourmelidis, Athanasios},
    title = {Discrete {$\Omega$}-results for the {R}iemann zeta function},
   journal = {Forum Math.},
    volume={37},
  number={4},
  pages={1221--1232},
  year={2025}
}

@article{soundararajan2008extreme,
  author = {K. Soundararajan},
  title = {Extreme values of zeta and {L}-functions},
  journal = {Math. Ann.},
  volume = {342},
  number = {2},
  pages = {467--486},
  year = {2008}
}

@article{qiyu2024large,
    title = {Large values of $\zeta(s)$ for $1/2<${R}e$(s)<1$},
journal = {J. Number Theory},
volume = {254},
pages = {199-213},
year = {2024},
author ={Qiyu Yang}
}

@article{XiaoYang2022,
  author = {X. Xiao and Q. Yang},
  title = {A note on large values of ${L}(\sigma, \chi)$},
  journal = {Bull. Aust. Math. Soc.},
  volume = {105},
  number = {3},
  pages = {412–418},
  year = {2022}
}

@article{yang2023omega,
  title={Omega theorems for logarithmic derivatives of zeta and {L}-functions},
  author = {D. Yang},
  journal={Preprint, arXiv:2311.16371},
  year={2023}
}

@article{yang2022logtype,
    title={ A note on log-type {GCD} sums and derivatives of the {R}iemann zeta function},
     author={Daodao Yang},
  journal={Preprint arXiv:2201.12968},
    year = {2022}
}
\end{document}